\documentclass{amsart}
\usepackage{graphicx}
\newtheorem{thm}{Theorem}[section]

\newtheorem{lem}[thm]{Lemma}

\theoremstyle{definition}
\newtheorem{defn}[thm]{Definition}
\theoremstyle{remark}
\newtheorem{rem}[thm]{Remark}
\numberwithin{equation}{section}
\begin{document}

\title{A Chevalley's theorem in class $ \mathcal{C}^r$}%

\author{G\'erard P. Barban\c con}%

\address{Department of Mathematics \\
    University of Texas at Austin \\
    1, University Station C1200\\
    Austin, TX 78712, USA\\}%

\email{gbarbans@math.utexas.edu}%

\thanks{I am grateful to Guy Laville (Univ. of Caen) and Hubert Rubenthaler
(Louis Pasteur Univ.) who invited me to give a talk at their
seminars from which I obtained a really useful feedback. I thank the
Univ. of Texas at Austin where I enjoy a most stimulating
environment and a wonderful library.}%

\subjclass{57R45, 13A50, 20F55}%
\keywords{Finite reflections groups, Coxeter groups, Whitney
 functions of class $\mathcal{C}^r$, formally holomorphic Whitney fields,
 Whitney extension theorem.}%

%\date{}%
%\dedicatory{}%
%\commby{}%
% ----------------------------------------------------------------
\begin{abstract}
 Let $W$ be a finite reflection group acting orthogonally on
 ${\bf R}^n$, $P$ be the Chevalley polynomial mapping determined
 by an integrity basis of the algebra of $W$-invariant polynomials,
 and $h$ be the highest degree of the coordinate polynomials in $P$.
 There exists a linear mapping:
 $\mathcal{C}^r(\mathbf{R}^n)^W \ni f \mapsto F\in \mathcal{C}^{[r/h]}(\mathbf{
 R}^n)$ such that $f=F\circ P$, continuous for the natural Fr\'echet
 topologies. A general counterexample shows that this result is
 the best possible. The proof uses techniques
 of division by linear forms and a study of
 compensation phenomenons. An extension to $P^{-1}(\mathbf{R}^n)$
 of invariant formally holomorphic regular fields is needed.
\end{abstract}
%----------------------------------------------------------------
\maketitle
% ----------------------------------------------------------------
\section{Introduction}
Let $W$ be a finite subgroup of $ O(n)$ generated by reflections.
The algebra of $W$-invariant polynomials is generated by $n$
algebraically independent $W$-invariant homogeneous polynomials and
the degrees of these basic invariants are uniquely determined. This
theorem and its converse were first stated by Shephard and Todd in
\cite{20} where the direct statement was proved case by case. Soon
after, in \cite{6} Chevalley gave a beautiful unique proof of this
direct statement, which is often called `Chevalley's theorem'.

A $W$-invariant complex analytic function may be written as a
complex analytic function of the basic invariants (\cite{21}).
Glaeser's theorem (\cite{11}) shows that real $W$-invariant
functions of class ${\mathcal C}^{\infty}$, may be expressed as
${\mathcal C}^{\infty}$ functions of the basic invariants.

In finite class of differentiability, Newton's theorem in class
${\mathcal C}^r$ (\cite{2}) dealt with symmetric functions and as a
consequence with the Weyl group of $A_n$. This particular case shows
a loss of differentiability as already did Whitney's even function
theorem (\cite{22}) which established the result for the Weyl group
of $A_1$. A first attempt to study the general case may be found in
the first part of \cite{4} where the best result was obtained for
the Weyl groups of $A_n, B_n$ and the dihedral groups $I_2(k)$ by a
method which was on the right track but needed an additional
ingredient to deal with the general case. The behavior of the
partial derivatives of the functions of the basic invariants on the
critical image was studied in \cite{12} for $A_n,B_n, D_n$ and
$I_2(k)$.

Here we give for any reflection group a result which is the best
possible as shown by a general counter example. Let $p_1, \ldots ,
p_n$ be an integrity basis, we define the `Chevalley' mapping
$P:{\mathbf R}^n \ni x\mapsto P(x)=(p_1(x), \ldots , p_n(x))\in
{\mathbf R}^n$. The loss of differentiability is governed by the
highest degree of the basic invariant polynomials. More precisely we
have:

\begin{thm}\label{thm1}
Let $W$ be a finite group generated by reflections acting
orthogonally on ${\mathbf R}^n$ and let $f$ be a $W$-invariant
function of class ${\mathcal C}^r$ on ${\mathbf R}^n$. There exists
a function $F$ of class ${\mathcal C}^{[r/h]}$ on ${\mathbf R}^n$
such that $f=F\circ P$, where $P$ is a Chevalley polynomial mapping
associated with $W$ and $h$ is the highest degree of the coordinate
polynomials in $P$, equal to the greatest Coxeter number of the
irreducible components of $W$.
\end{thm}

  Since a change of basic invariants is an invertible polynomial map
  on the target, the statement does not depend on the choice of the
  set of basic invariants.
\section{ The Chevalley mapping}

A detailed study may be found in \cite{5} or \cite{8}.

  When $W$ is reducible, it is a direct product of its irreducible
 components, say $W= W^0\times W^1 \times \ldots \times W^s $ and we may
 write $\mathbf{R}^n$ as an orthogonal direct sum $\mathbf{R}^{n_0}\oplus
 \mathbf{R}^{n_1}\oplus \ldots \oplus \mathbf{R}^{n_s}$ where $W^0$ is
 the identity on $\mathbf{R}^{n_0}$,
 subspace of $W$-invariant vectors, and for $i=1,\ldots, s$, $W^i$
 is an irreducible finite Coxeter group acting on ${\bf R}^{n_i}$.
 We may choose coordinates that fit with this orthogonal direct sum.
 If $w=w_1\ldots w_s\in W$ with $w_i\in W^i, \; 1\le i\le s$
 we have $w(x)= w(x^0, x^1,\ldots , x^s)= (x^0, w_1(x^1),\ldots,
 w_s(x^s))$ for all $x\in \mathbf{R}^{n}$. The
 direct product of the identity on ${\bf R}^{n_0}$ and Chevalley
 mappings $P^i$ associated with $W^i$ acting on $\mathbf {R}^{n_i},\; 1
 \le i\le s$, is a Chevalley map $P= Id_0\times P^1\times \ldots
 \times P^s$ associated with the action of $W$ on $\mathbf{R}^{n}$.

  For an irreducible $W$ (or for an irreducible component)
 we will assume as we may that the degrees of the coordinate
 polynomials $p_1,\ldots ,p_n$ are in increasing order:
 $2=k_1\le \ldots \le k_n=h$, Coxeter number of $W$.
 In the reducible case we will denote by $h$ the highest degree of
 the coordinate polynomials which is the maximal Coxeter number of
 the irreducible components.

Let $\mathcal{R}$ be the set of reflections different from identity
in $W$. The number of these reflections is
$\mathcal{R}^{\#}=d=\sum_{i=1}^n(k_i-1)$. For each $\tau \in
\mathcal{R}$, let $ \lambda_{\tau}$ be a linear form the kernel of
which is the hyperplane $ H_{\tau}=\{ x\in {\bf R}^n |\tau (x) = x
\} $. The jacobian of $P$ is $J_P=c \prod_{\tau \in  \mathcal{R}}
\lambda_{\tau}$ for some constant $c\neq 0$. The critical set is the
union of the $ H_{\tau}$ when $\tau$ runs through $ \mathcal{R}$.

 A Weyl Chamber $C$ is a connected component of the regular
 set. The other connected components are obtained by the
 action of $W$ and the regular set is $\bigcup_{w\in W} w(C)$.
 There is a stratification of ${\mathbf R}^n$ by the regular set, the
 reflecting hyperplanes $ H_{\tau}$ and their intersections.
 The mapping $P$ induces an analytic diffeomorphism of $C$ onto
 the interior of $P({\mathbf R}^n)$ and an homeomorphism
 that carries the stratification from the
 fundamental domain $\overline{C}$ onto $P(\mathbf{R}^n)$.

 The Chevalley mapping is neither injective nor surjective. Actually
 the fiber over each point of the image is a $W$-orbit. The mapping $P$
 is proper and separates the orbits (\cite{19}). It is the restriction
 to $\mathbf{R}^n$ of a complex $W$-invariant mapping from $\mathbf{C}^n$
 onto (\cite{13}) $\mathbf{C}^n$, still denoted by $P$.

 \noindent On its regular set, the complex $P$ is a local analytic isomorphism.
 Its critical set, where the jacobian vanishes, is the union of the
 complex hyperplanes $H_{\tau} =\{ z\in {\mathbf C}^n | \tau(z)=z \}$,
 kernels of the complex forms $\lambda_{\tau}$.
 The critical image is the algebraic set $\{ u\in {\mathbf C}^n |
 \Delta(u)= J^2_P(z)=0 \}$, on which $P$ carries the stratification.
 The complex $P$ is proper, it is a $W^{\#}$-fold covering of ${\bf C}^n$
 ramified over the critical image.

 Finally, there are only finitely many types of irreducible finite Coxeter groups
 defined by their connected graph types. Even when these groups are Weyl
 groups of root systems or Lie algebras, we will follow the general usage
 and denote them with upper case letters.

 \section{ Whitney Functions and $r$-regular, $m$-continuous jets.}

A complete study of Whitney functions may be found in
\cite{21}.\vskip 5pt

A jet of order $m\in {\mathbf N}$, on a locally closed set $E\subset
{\mathbf R}^n$ is a collection $A=(a_k)_{k\in {\bf N^n}\atop \mid
k\mid \le m}$ of real valued functions $a_k$ continuous on $E$. At
each point $x\in E$ the jet $A$ determines a polynomial $A_x(X)$,
and we sometimes speak of continuous polynomial fields instead of
jets (\cite{15}). As a function, $A_x$ acts upon vectors $x'-x$ tangent
to ${\bf R}^n$ at $x$. To avoid introducing the notation $T_x^r A$,
we write somewhat inconsistently:
\[ A_x: x' \mapsto A_x(x')=\sum_k {1\over
k!} a_k(x)\; (x'-x)^k.\]

\noindent By formal derivation of $A$ of order $q\in {\mathbf
N}^n,\; \mid q\mid \le m$, we get jets of the form $(a_{q+k})_{\mid
k\mid \le m-\mid q\mid}$ inducing polynomials \[(D^qA)_x(x')= \left(
{\partial^{\mid q\mid} A\over
\partial x^q}\right)_x (x')=a_q(x)+\sum_{k>q\atop |k|\le m}{1\over (k-q)!}
a_k(x)\; (x'-x)^{k-q}.\]

\noindent For $0\le \mid q\mid \leq r\le m$, we put:
$$(R_xA)^q(x')=(D^q A)_{x'}(x')-(D^q A)_x(x').$$

\begin{defn} Let $A$ be an $m$-jet on $E$. For $r\le m , A$ is
$r$-regular on $E$, if and only if for all compact set $K$ in $E$,
for $(x, x')\in K^2$, and for all $q\in {\bf N}^n$ with $\mid q\mid
\leq r$, it satisfies the Whitney conditions.

\centerline {$({\mathcal W}_q^r) \qquad (R_xA)^q(x') = o (\mid
x'-x\mid^{r-\mid q\mid})$, when $\mid x-x'\mid \to 0$.}
\end{defn}

 \begin{rem} Even if $m> r$ there is no need to
consider the truncated field $A^r$ in stead of $A$ in the conditions
$({\mathcal W}_q^r)$. Actually $(R_xA^r)^q(x')$ and $(R_xA)^q(x')$
differ by a sum of terms $[a_k(x)/(k-q)!] \, (x'-x)^{k-q}$, with
$a_k$ uniformly continuous on $K$ and $ \vert k\vert -\vert q\vert
> r-\vert q\vert $.\end{rem}

The space of $r$-regular jets of order $m$ on $E$, is naturally
provided with the Fr\'echet topology defined by the family of
semi-norms: \[\Vert A\Vert_{K_n}^{r,m}=\sup_{x\in K_n \atop \mid
k\mid \leq m} \mid {1\over k!}a_k(x)\mid +\sup_{(x,x')\in K_n^2
\atop x\neq x', \mid k\mid \leq r} \left( {\mid (R_xA)^k(x')\mid
\over \mid x-x'\mid^{r-\mid k\mid}}\right)\] where $K_n$ runs
through a countable exhaustive collection of compact sets of $E$.
Provided with this topology the space of $r$-regular, $m$-continuous
polynomial fields on $E$ is a Fr\'echet space denoted by ${\mathcal
E}^{r,m}(E)$. If $r=m,\; {\mathcal E}^r(E)$ is the space of Whitney
fields of order $r$ or Whitney functions of class ${\mathcal C}^r$
on $E$.

\begin{thm} \label{Wextsion} {\rm Whitney extension theorem (\cite{23}}).
The restriction mapping of the space ${\mathcal E}^r({\mathbf R}^n)$
of functions of class ${\mathcal C}^r$ on ${\mathbf R}^n $ to the
space ${\mathcal E}^r(E)$ of Whitney fields of order $r$ on $E$, is
surjective. There is a linear section, continuous when the spaces
are provided with their natural Fr\'echet topologies. \end{thm}

Let $E$ be a closed subset of ${\mathbf C}^n\simeq {\mathbf
R}^{2n}$, we consider jets $A$ on $E$ with complex valued
coefficients $a_k$. They induce in $z\in E$ the polynomials:
\[A_z(X,Y)= \sum_{\mid k \mid +\mid l\mid \leq m} {1\over {k!l!}}
a_{k,l}(z)\; X^kY^l \in {\bf C}[X,Y],\] and we may define the
Fr\'echet space ${\mathcal E}^{r}(E;{\mathbf C})$ of complex valued
Whitney functions of class ${\mathcal C}^r$.

\begin{defn}(\cite{15})  A Whitney function $A\in {\mathcal E}^{r}(E;{\bf
C})$ is formally holomorphic if it satisfies the Cauchy-Riemann
equalities: \[i{\partial A\over\partial {X_j}}= {\partial A\over
\partial {Y_j}}, \; j=1, ... ,n.\] \end{defn}

 Let $Z=(Z_1,\ldots ,Z_n), \; Z_j=X_j+iY_j,\; j=1, .\ldots ,n$.
 The field $A$ is formally holomorphic if and only if
 $ \partial A/\partial \overline{Z_j} =0, \; j=1, ... ,n $.
 Thus for all $z\in E$ the polynomial $A_z$ belongs to ${\bf C}[Z]$
 and is of the form $\displaystyle{A_z(Z)=\sum_k {1\over k!} a_k(z)
 Z^k} $.\vskip 5pt

The algebra of formally holomorphic Whitney functions of class
${\mathcal C}^r$ on the (locally) closed set $E$ of ${\mathbf C}^n$
will be denoted by ${\mathcal H}^{r}(E)$. It is a closed sub-algebra
of ${\mathcal E}^{r}(E;{\mathbf C})$ and therefore a Fr\'echet space
when provided with the induced topology. In practice we define the
semi-norms $\Vert A\Vert^{K_n}_{r}$ on ${\mathcal H}^{r}(E)$ by the
same formulas as in ${\mathcal E}^{r}(E;{\mathbf R})$, only using
modulus instead of absolute value.

To take advantage of compensation phenomenons, it may be convenient
to consider Fr\'echet spaces ${\mathcal H}^{r,m}(E)$ of formally
holomorphic $r$-regular jets of order $m \ge r$ on $E$.

\begin{defn} A real form (\cite{17}) or real
 situated subspace (\cite{15}) of ${\bf C}^n$ is a real
 vector subspace $E$ of real dimension $n$ such that $E \oplus i E
 ={\bf C}^n$.\end{defn}

 \noindent A real form is a real subspace
 $E_S= \{z\in {\mathbf C}^n| Sz=z \}$, where $S$ is an anti-involution.

Example. Let $\alpha$ be an involution of $\{1,\ldots,n\}$,
 $\Gamma_{\alpha}=\{ z\in {\mathbf C}^n | z_{\alpha (j)}=
 \overline{z_j},\; j=1,\ldots,n\}$ is a real form of ${\mathbf C}^n$ defined by
 the anti-involution $z\mapsto \overline{\alpha( z)}$.

 Let $W$ be a finite reflection group acting orthogonally on ${\bf
 R}^n$ and $P$ be its Chevalley polynomial mapping as above.
 Since $P$ is defined over ${\mathbf R}$ (its coefficients are real),
 $P^{-1}({\mathbf R}^n)$ is the union of real forms $\Gamma_{S_w}
 \subset{\mathbf C}^n$, where $w$ runs through the involutions of $W$
 and $S_w$ is the anti-involution defined by $S_w(u+iv)= wu-iwv$.
 \vskip 5pt

 Let $f\in {\mathcal C}^r({\mathbf R}^n)^W$ be a $W$-invariant function of class
 ${\mathcal C}^r$. It induces on ${\bf R}^n$ a $W$-invariant Whitney
 field of order $r$ and
 a formally holomorphic field in ${\mathcal H}^{r}({\mathbf R}^n)^W$ which
 will still be denoted by $f$. By using Whitney's extension theorem,
 one may show (\cite{2}) that there is a linear and continuous extension
 \footnote{This extension will allow us to get a field $\tilde F$ on ${\bf R}^n$
 with $\tilde f=\tilde F\circ P$ and derive its regularity from its continuity
 by using the below lemma \ref{lem3}. This process might be avoided if there was an
 available proof of the Whitney regularity property (\cite {24}) of $P({\bf R}^n)$,
 a most likely conjecture proved for $A_n$ in \cite {14}, easy to show for some
 lower dimensional reflection groups (\emph{e.g.} $I_2(k), H_3$). Unfortunately there is none, although in \cite{9} a key point for the proof of this conjecture is studied and a sketchy proof of it is given for several Coxeter groups.}:
 \[{\mathcal H}^r({\mathbf R}^n)^W \ni f \mapsto \tilde f \in {\mathcal H}^r
 (P^{-1}({\mathbf R}^n))^W .\]

\section{Some multiplication and division properties.}

\begin{lem}\label{lem1} Let $\Gamma$ be a finite union of real forms of ${\bf
C}^n$, $A$ be in $\mathcal{H}^{r}(\Gamma)$, and $Q$ be a polynomial
$(s-1)$-flat on $S$. Let $z\in \Gamma$ and $z_0\in S\cap \Gamma$,
then for all $q\in {\bf N}^n, \mid q\mid \le r$:
\[(R_{z_0}QA)^q(z)=(D^q QA)_{z}(z)-(D^q QA)_{z_0}(z)\in o(\mid
z-z_0\mid^{r-\mid q\mid +s}).\]  Moreover $QA\in \mathcal
{H}^{r+s}(S\cap \Gamma)$ and is $(s-1)$-flat on $S\cap \Gamma$. For
all compact $K\subset S\cap\Gamma$,
 there exists a numerical constant $c$ such that $\|QA\|^{r+s}_K \le
 c\|Q\|^{r+s}_K\|A\|^r_K$. \end{lem}
{\it Proof.}  Let $z_0\in S\cap \Gamma$. For all $z\in \Gamma$, all
$q\in {\bf N}^n, \mid q\mid \le r$, and $p\le q$, we consider:
\[(D^p Q)_{z}(z)(D^{q-p} A)_{z}(z)-(D^{ p} Q)_{z_0}(z)(D^{q-p} A
)_{z_0}(z).\] By Taylor's formula for polynomials $( D^{p}
Q)_{z}(z)=(D^{p} Q)_{z_0}(z)$, and this difference is: \[(D^{p}
Q)_z(z)\left[(D^{q-p} A)_{z}(z)-(D^{ q-p} A)_{z_0}(z) \right].\] By
assumption $(D^{q-p} A)_{z}(z)-( D^{q-p} A )_{z_0}(z)\in o( \mid
z-z_0 \mid^{r-\mid q\mid +\mid p \mid } )$, and for $\mid p\mid  <
s$ $(D^{p} Q)_z(z)\in O(\mid z-z_0 \mid^{s- \mid p\mid })$. The
product is in $o(\mid z-z_0 \mid^{ r-\mid q\mid +s})$ either because
$|p|<s$ and $r-|q|+|p|+s-|p|=r-|q|+s$ or because $|p|\geq s$ and
$r-|q|+|p|\geq r-|q|+s$.\vskip 5pt

\noindent The behavior of $(R_{z_0}QA)^q(z)$ is now a consequence of
Leibniz' derivation formula.\vskip 5pt

Actually $QA\in \mathcal{H}^{r,r+s}$. On $S\cap \Gamma$, $\mid p\mid
< s \Rightarrow (D^{p} Q)_{z_0}(z_0)=0$, therefore in the
derivatives of $QA$ of order $\le r+s$ the only derivatives of $A$
that are not multiplied by a derivative of $Q$ that vanishes, are of
order $\le r$. Then the above estimates show that when $|q|\le r+s$,
the field $QA$ satisfies Whitney conditions $\mathcal{W}_q^{r+s}$ on
$S\cap \Gamma$. Thus $QA\in \mathcal {H}^{r+s}(S\cap \Gamma)$, and
clearly it is $(s-1)$-flat on $S\cap \Gamma$.

This situation was already noticed in \cite{10}: when multiplying a field
$r_1$-regular, $(s_1-1)$-flat by a field $r_2$-regular,
$(s_2-1)$-flat, the product is $ \min (r_1+s_2, r_2+s_1)$-regular
and $(s_1+s_2-1)$ flat. Here, on $S\cap \Gamma$, we have $r_1=r,\;
s_1=0$ for $A$ and $\;r_2= +\infty,\; s_2=s$ for $Q$. $\;\Box$
\vskip 10pt

Example.  Let $Q$ be an homogeneous polynomial of degree $s$. It
vanishes at the origin with all its derivatives of order $\le s-1$.
If $A\in \mathcal{H}^{r}(\Gamma)$, for all $z\in \Gamma$ and all
$q\in {\bf N}^n, \mid q\mid \le r$: \[(R_{0}QA)^q(z)=(D^{ q}
QA)_{z}(z)-(D^{q} QA)_{0}(z)\in o(\mid z\mid^{r+s-\mid q\mid}).\]

\noindent The same result holds if instead of a product $QA$ we have
a sum $\sum_{i=1}^{n}Q_iA_i$, with homogeneous polynomials $Q_i$ of
degree $s_i\geq s$ and $A_i\in \mathcal{ H}^{r}(\Gamma)$.

Let us recall the following division lemma:

\begin{lem} [\cite{2}] \label{lem2} Let $\Gamma$ be a finite union of real forms of
${\bf C}^n$, and $\lambda\neq 0$ be a complex linear form with
kernel $H$. If $A\in \mathcal{H}^{r}(\Gamma)$ is such that $A_z(Z)$
is divisible by $\lambda_z(Z)$ whenever $z\in \Gamma \cap H$ then
there exists a field $B\in \mathcal{H}^{r-1}(\Gamma)$ such that $A^r
= (\lambda B )^r$. For all compact $K\subset \Gamma$, there exists a
constant $c$ such that $\|B\|^{r-1}_K \le c\|A\|^r_K$.\end{lem}

\noindent Actually $B\in \mathcal{H}^{r}(\Gamma\setminus H )$ and if
$\mid s\mid = r$, then $\displaystyle{\lambda (z)(D^{s} B)_{z}(z)}$
tends to zero with $\lambda (z)$.\vskip 5pt

\begin{rem} The lemma still holds if we replace $\Gamma$ by its
intersection with one or several hyperplanes distinct from $H$.
\end{rem}

The proof of lemma \ref{lem2} relies upon a consequence of the mean
value theorem that will be instrumental in what follows: \vskip 5pt

 \begin{lem}[\cite{1},\cite{15}]\label{lem3} Let $\Gamma$ be a finite union of
 real forms of ${\bf C}^n$, $\Delta\ne 0$ be a polynomial,
 and $X= \{x\in {\bf C}^n  \mid \Delta(x)=0 \}$.
 If $f\in \mathcal{H}^{r}(\Gamma \setminus X)$ is $r$-continuous on
 $\Gamma$, then $f\in \mathcal{H}^{r}(\Gamma)$. \end{lem}

 Let $(\lambda_{\tau})_{\tau\in \mathcal{D}}$ be
 $\mathcal{D}^{\#}=d$ non zero complex linear forms with kernels
 $(H_{\tau})_{\tau\in \mathcal{D}}$ all distinct. The hyperplanes
 $(H_{\tau})_{\tau\in \mathcal{D}}$ and their intersections induce a
 stratification on $\Gamma$. Let $S_p$ be a stratum, connected component
 of the intersection of $\Gamma$ and exactly $p$ of these hyperplanes, say
 $ (H_{\tau})_{\tau\in \mathcal{B}_p}$, $\mathcal{B}_p^{\#}=p$.
 The border $\overline{S_p}\setminus S_p$ is a union
 $\bigcup S_{p+l}$ of strata of lower
 dimensions, containing $S_d=\Gamma\cap (\bigcap_{\tau\in \mathcal{
 D}}H_{\tau})$. Using these notations we have:\vskip 5pt

 \begin{lem}\label{lem4} For $i=1,\ldots,n$, let $A_i$ be in $\mathcal{H}^r(\Gamma)$ and $Q_i$
 be an homogeneous polynomial $(s_p-1)$-flat on
 $S_p$ and more generally $(s_{p+l}-1)$-flat on each of the
 $S_{p+l}$. Assume $p+l-s_{p+l}$ is an increasing function of $l$
 and that $A=\sum_{i=1}^n Q_iA_i=(\prod_{\tau \in \mathcal{D}}
 \lambda_{\tau}) C$,  meaning that:
 \[\forall \mathcal{U}\subseteq \mathcal{D}, \; \; A_z(Z) {\it is\; divisible\; by}
 \prod_{\tau \in \mathcal{U} }\lambda_{\tau}(Z)\;{\it when}\; z\in \Gamma \cap
 (\bigcap_{\tau \in \mathcal{U}} H_{\tau}).\] The field $C$
 is in $\mathcal{H}^{r+s_p-p} (S_p)$ and its coefficients of order $\le r+s_d-d$
 are continuous on $ \overline{S_p}$. \footnote{Actually since $\overline{S_p}$ is
 convex and thus Whitney 1-regular (\cite{24}), lemma \ref{lem4} yields that $C$
 is $(r+s_d-d)$-regular on $\overline{S_p}$.}
 \end{lem}

 {\it Proof.} By lemma \ref{lem1}, $\sum_{i=1}^n Q_iA_i$ is in ${\mathcal H}^{r+s_p}(S_p)$
  and in ${\mathcal H}^{r+s_{p+l}}(S_{p+l})$.
  By lemma \ref{lem2}, the field $C$ is in
  ${\mathcal H}^{r+s_p-p}(S_p)$ and in ${\mathcal H}^{r+s_{p+l}-(p+l)}(S_{p+l})$.
  We are just to show the continuity on $\overline{S_p}$ of the coefficients
  of order $\le r+s_d-d$ in $C$.

  Let $S_{p+q}$ be one of the strata of largest dimension in $\overline{S_p}
  \setminus S_p$, and let ${\mathcal B}_q$ with ${\mathcal B}_q^{\#}=q$,
  be the subset of ${\mathcal D}$
  such that $S_{p+q}$ is a connected component of the intersection of
  $\Gamma$ and the hyperplanes $(H_{\tau})_{\tau\in {\mathcal B}_p\cup{\mathcal B}_q}$,
  but no other. We may have $q=1$ but not necessarily since the addition of one hyperplane
  may automatically entail the addition of some more.  \vskip 5pt

   We put:
  $A= (\prod_{\tau\in {\mathcal B}_p} \lambda_\tau)( \prod_{\tau\in {\mathcal B}_q}
  \lambda_\tau)( \prod_{\tau\in {\mathcal D}\setminus({\mathcal B}_p\cup{\mathcal B}_q)}
  \lambda_\tau) \; C$
  and define:

  \centerline{$C^1=(\prod_{\tau\in {\mathcal D}\setminus({\mathcal B}_p\cup{\mathcal B}_q)}
  \lambda_\tau) \; C ,\quad {\rm and}
  \quad B=(\prod_{\tau\in {\mathcal B}_q} \lambda_\tau) \; C^1.$}

 \noindent On $S_p$, $B$ is in ${\mathcal H}^{r+s_p-p}$ and so are $C$ and
 $C^1$. On $S_{p+q}$, $C$ and $C^1$ are in ${\mathcal H}^{r+s_{p+q}-(p+q)}$.

 \noindent Let $z_0$ be the orthogonal projection on $S_{p+q}$ of some $z\in
 \Gamma$. By lemma \ref{lem1}, we have:
 \[A_z(z)-A_{z_0}(z)=[\prod_{\tau \in {\mathcal B}_p}\lambda_{\tau}(z)]
  [B_z(z)-B_{z_0}(z)] \in o(\mid z-z_0\mid^{r+s_{p+q}}).\]
  Let $\pi$ be a derivation of order $p$, by Leibniz derivation formula:
   \[D^{\pi}A_z(z)-D^{\pi}A_{z_0}(z)=[\prod_{\tau \in {\mathcal
   B}_p}\lambda_{\tau}(z)]
  [D^{\pi}B_z(z)-D^{\pi}B_{z_0}(z)]+\ldots \] \[\ldots+
  c \;[B_z(z)-B_{z_0}(z)]\in o(\mid z-z_0\mid^{r+s_{p+q}-p})\] for
  some constant $c\neq 0$.

  For the remaining part of the proof, we assume that $z$ is in $S_p$.
  Then for all $\tau\in {\mathcal B}_p,\; \lambda_\tau (z)=0$ and we get:
  \[B_z(z)-B_{z_0}(z) \in o(\mid z-z_0\mid^{r+s_{p+q}-p}).\]

 More generally, by considering derivations of order $\pi +\kappa$,
 with $\mid \kappa \mid =k \le r+s_p-p$ \;we would get in the same way:
 \[D^{\kappa}B_z(z)-D^{\kappa}B_{z_0}(z)\in o(\mid
 z-z_0)\mid^{r+s_{p+q}-p-k}).\]

Let us put $E=\Gamma\cap (\bigcap_{\tau\in {\mathcal B}_p}
H_{\tau})$ and $F = H_{\tau}$ for some $\tau \in {\mathcal B}_q$.
Since $E$ and $F$ are ${\bf R}$-linear subspaces, there exists a
constant $a> 0$ such that for all $ x ,\; d(x,E)+d(x,F) \ge a\; d(x,
E\cap F)$ or equivalently, a constant $b> 0$ such that for all $
x\in E$, $d(x,F) \ge b\; d(x, E \cap F)$. We say that $E$ and $F$
are regularly separated. \footnote{This is a trivial case of regular
separation and does not need a general study. The interested reader
may take a look at \cite{16} or \cite{21}. The regular separation of
real forms was implicitly used in the above extension of $f$ to
$P^{-1}(\mathbf{R}^n)$}

\noindent The regular separation brings the existence of constants
$c_{\tau}$ such that: \[\forall\tau\in {\mathcal B}_q,\quad
|z-z_0|\le c_{\tau}\; d(z,H_{\tau})=
c_{\tau}^1|\lambda_{\tau}(z)|.\] Therefore, since
\[B_z(z)-B_{z_0}(z)= [\prod_{\tau \in {\mathcal
B}_q}\lambda_{\tau}(z)] [C^1_z(z)-C^1_{z_0}(z)]\in o(|
z-z_0|^{r+s_{p+q}-p}),\] we get that $C^1_z(z)-C^1_{z_0}(z) \in
o(\mid z-z_0\mid^{r+s_{p+q}-(p+q)})$.\vskip 5pt

\noindent Let us assume by induction that for $\mid l\mid \;\le  k
-1\; <r+s_{p+q}-(p+q)$ : \[D^l C^1_z(z)-D^lC^1_{z_0}(z) \in o(\mid
z-z_0\mid^{r+s_{p+q}-(p+q)-\mid l\mid}).\] By Leibniz derivation
formula, for $j,\;|j|=k$: \[D^jB_z(z)-D^jB_{z_0}(z) =[\prod_{\tau
\in {\mathcal B}_q}\lambda_{\tau}(z)](D^j C^1_z(z)-D^j
C^1_{z_0}(z))+\]
 \[+ \sum_{k-q\le |j_i|= k-l \le k-1}
 a_{q-l}[\prod_{q-l}\lambda_{\tau}(z)]
(D^{j_i}C^1_z(z)-D^{j_i}C^1_{z_0}(z))\in o(\mid
z-z_0\mid^{r+s_{p+q}-p- k}).\]

\noindent The $[\prod_{q-l}\lambda_{\tau}(z)]$ stand for $
D^{j-j_i}(\prod_{\tau \in {\mathcal B}_q}\lambda_{\tau})(z)$, up to
a constant factor included in $a_{q-l}$. Applying the induction
assumption to the derivations $D^{j_i}$ of order $\le k -1$, we see
that each term of the sum is in \[o(\mid
z-z_0\mid^{r+s_{p+q}-(p+q)-( k-l)+q-l}) = o(\mid
z-z_0\mid^{r+s_{p+q}-p-k})\] and thus, that the first term also is.
Then, using the regular separation as above, we obtain: \[D^j
C^1_z(z)-D^j C^1_{z_0}(z)\in o(\mid z-z_0\mid^{r+s_{p+q}-(p+q)-
k}).\] This completes the induction, and shows that the coefficients
of $C^1$ of order $\le r+s_{p+q}-(p+q)$ are continuous in $z_0$.
Since they are continuous on $S_{p+q}$, by using the triangular
inequality, we get their continuity on $S_p \cup S_{p+q}$.

 \noindent Any $z_1\in S_{p+q}$ has a neighborhood which does not
 meet any of the $H_{\tau}$ but those containing $z_1$ and
 in this neighborhood, the continuity of the
  coefficients of $C$ and $C^1$ are the same. Hence the continuity on
  $S_p\cup S_{p+q}$ of the coefficients of $C$ of order $\le r+s_{p+q}-(p+q)$.

  \noindent We can get an analogous result for each
  stratum $S_{p+q'}$ of maximal dimension in
  $\overline{S_p}\setminus S_p$ ($q $ and $q'$ may be equal but ${\mathcal B}_q $ and
  ${\mathcal B}_{q'}$ are distinct).

  \noindent We then proceed with the strata of highest dimension in
  $\overline{S_{p+q}}-S_{p+q}$ to get the continuity on $S_{p+q}\cup
  S_{(p+q)+l}$ of the coefficients of order at most $r+s_{p+q+l}-(p+q+l)$.
  The continuity on $S_p\cup S_{p+q}\cup
  S_{(p+q)+l}$ when $z\in S_p$ tends to $z_1\in S_{(p+q)+l}$ is obtained by the
  triangular inequality, considering the orthogonal projection $z_0$
  of $z$ on $S_{p+q}$. We can do the same for each of the
  $\overline{S_{p+q'}}-S_{p+q'}$, and continue
  the processes until we reach $S_d$ through all possible paths, getting the
  continuity of coefficients of order at most $r+s_d-d$. Thus we
  have the continuity of these coefficients between any two points
  in $S_{p+l}$ and $S_{(p+l)+m}$.

  \noindent About the continuity between points in $\overline{S_{p+q}}$ and
  $\overline{S_{p+q'}}$, if they are in the intersection we already have
  the continuity, but if they are not, we consider their orthogonal
  projections on the intersection and thanks to the regular separation of
  the strata, we can get the continuity by the triangular inequality.

 \noindent Then the global $(r+s_d-d)$-continuity on $\overline{S_p}$ is
 clear. $\quad \Box$

  \begin{rem} When $p=0$, the strata of type $S_0$ are open
  and $s_0=0$. For $q=1$, the first step is given by
  lemma \ref{lem2}.\end{rem}

\section{ Proof of Theorem 1.1.}

  We consider an invariant function
  $ f\in \mathcal{C}^r({\bf R}^n)^W$. This function or rather the formally
  holomorphic field it induces on ${\bf R}^n$ has a linear and
  continuous extension $\tilde f \in \mathcal{H}^r(P^{-1}({\bf R}^n))^W$.

 \subsection*{Pointwise solution.}

\begin{lem} [\cite{4}]\label{lem5} For all $W$-invariant, formally
 holomorphic polynomial field $\tilde f$ of degree $r$ on $P^{-1}({\bf
 R}^n)$, there exists a formally holomorphic polynomial field
$\widetilde F$ of degree at most $r$ such that for all $z\in P^{-1}
 ({\bf R}^n),\; \tilde f_z= (\widetilde F_{P(z)}\circ P)^r_z $.
 \end{lem}

 {\it Proof.} On the complement of
$\Gamma \cap \bigcup_{\tau\in \mathcal{R}} H_{\tau}$ in $\Gamma $,
the mapping $P$ is a local analytic isomorphism and this yields the
construction of $\widetilde F= (\tilde f\circ P^{-1})^r$,
unambiguously since both $ \tilde f$ and $ P$ are $ W$-invariant. On
the regular image of $P$, $\widetilde F$ verifies $\tilde
f^r=(\widetilde F\circ P)^r$.

Let $x\in \Gamma \cap (\bigcup_{\tau\in \mathcal{R}}H_{\tau})$ and
let $W_x$ be the isotropy subgroup of $W$ at $x$. The polynomial
$\tilde f_x$ is $W_x$-invariant since for all $w_0\in W_x\subset W$:
$ \tilde f_x(X)=\tilde f_{w_0x}(w_0X)=\tilde f_x(w_0X)$ where the
first equality results from the $W$-invariance of the field $\tilde
f$ and the second from $w_0x=x$.
 As a consequence, $\tilde f_x$ is a polynomial in the
$W_x$-invariant generators $v=(v_1, \ldots , v_n)$ of the subalgebra
of $W_x$-invariant polynomials, and we have $\tilde f_x=Q\circ v$.

There exists a neighborhood of $x$ in ${\bf C}^n$ which does not
meet any of the hyperplanes $H_{\tau}$ but those containing $x$. In
this neighborhood we may write $P= q\circ v$ for some polynomial
$q$, since $P$ is $W_x$-invariant. Up to a multiplicative constant
the jacobian of $q$ at $ v(x)$ is the product
$\prod_{\lambda_s(x)\neq 0} \lambda_s$ and $q$ is an analytic
isomorphism in a neighborhood of $v(x)$.

We define the jet at $P(x)$ by $\widetilde F_{P(x)}= [Q\circ
q^{-1}]^r$ and get:

$[\widetilde F\circ P]^r_x= [(Q\circ q^{-1})^r\circ (q\circ v)]^r_x=
[(Q\circ q^{-1})\circ (q\circ v)]^r_x=(Q\circ v)_x= \tilde f_x .
\quad \Box$

\begin{rem} When the isotropy subgroup of $x_0$ is $W$ itself,
$\; \forall w\in W, \; \tilde f_{x_0}(X)=\tilde f_{wx_0}(wX)=\tilde
f_{x_0}(wX)$. This means that $\tilde f_{x_0}(X)$ is $W$-invariant
and by the polynomial Chevalley's theorem, that $\tilde
f_{x_0}(X)=Q_0(P(X))$. The polynomial $Q_0=\widetilde F_{P(x_0)}$ of
weight $r$ is of degree $[r/h]$ in the invariant polynomial $p$ of
highest degree $h$. The result announced in theorem \ref{thm1} fits
with the formal computation.\end{rem}

\subsection*{A criterion of regularity for $\tilde{F}$}

 When $\tilde f\in \mathcal{H}^r(P^{-1}(\mathbf{R}^n))$, the above
 proof of lemma \ref{lem5} shows that $\widetilde F= (\tilde f\circ P^{-1})^r$ is
 r-regular on the complement in ${\bf R}^n$ of
 the critical image $\{u\in {\bf C}^n\mid \Delta (u)=0\}$. The
 discriminant $\Delta $ is a polynomial, therefore by lemma \ref{lem4}
 it will be sufficient to prove that $\widetilde F$ is $[r/h]$-continuous
 on ${\bf R}^n$ to get its $[r/h]$-regularity.

 Since $P$ is proper the continuity of $\widetilde F_{\alpha}\circ P$,
 entails the continuity of the coefficient $\widetilde F_{\alpha}$.
 So let us check the continuity of the $\widetilde F_{\alpha}\circ P$
 when $\mid \alpha \mid \le [r/h]$.

 Clearly $\widetilde F_0\circ P=\tilde f_0$ is continuous.
 For the first derivatives, it is natural to consider the partial
 derivatives of $\tilde f$, and get the system:
 \begin{equation}
  \left( {\partial \tilde f\over \partial z}\right) = \left(\left(
 {\partial p_i\over \partial z_j}\right)_{1\le i\le n\atop 1\le j\le n}\right)
 \left({\partial \widetilde F\over \partial p} \circ P\right). \label{sys1}
 \tag{I}
 \end{equation}

 \noindent If we show that the loss of differentiability from
 $\tilde f=\widetilde F\circ P$ to $(\displaystyle{{\partial \widetilde
 F/ \partial p}) \circ P}$ when solving ~\eqref{sys1} is of $h$ units, applying
 the same process to $\tilde g_j=({\displaystyle{\partial \widetilde
 F}/ \displaystyle{\partial p_j}})\circ P$ instead of $\tilde f=\widetilde F \circ
 P$, at the next step again there will be a loss of differentiability
 of $h$ units. An induction would show that for
 $|\alpha|\le \displaystyle{[{r\over h}]}$, the mappings
 $({\displaystyle{\partial^{|\alpha|} \widetilde F}/
 \displaystyle{\partial p^{\alpha}}})\circ P$ are continuous on
 $P^{-1}({\bf R}^n)$ and since $P$ is proper, that the derivatives
 ${\displaystyle{\partial^{|\alpha|} \widetilde F}/
 \displaystyle{\partial p^{\alpha}}}$ of $\widetilde F$
  are continuous on ${\bf R}^n$.\vskip 5pt

 \noindent \emph{Conclusion:} To complete the proof, we just have to show that
 when solving ~\eqref{sys1}, $\forall j=1, \ldots,n$ we get:
 $({\displaystyle{\partial \widetilde F}/ \displaystyle{\partial
 p_j}})\circ P \in \mathcal{H}^{r-h}(P^{-1}({\bf R}^n))$.

\subsection*{Solving ~\eqref{sys1}, Reduction to the irreducible
case}

 \noindent Using Cramer's method, we multiply both sides by the adjoint matrix of the system.
 Since the jacobian determinant is
 $ c \;\displaystyle{(\prod_{\tau\in \mathcal{R}}\lambda_{\tau})}$, we get:
 \begin{equation}
 \left \{ c \;\displaystyle{(\prod_{\tau\in \mathcal{R}}
 \lambda_{\tau})\; \;{\partial \widetilde F\over
 \partial p_j}\circ P =\sum_{i=1}^n
 (-1)^{i+j}M_{i,j}{\partial \tilde f \over \partial
 z_i.}},\; j=1\ldots ,n \right. \label{soleq2}\tag{II}
 \end{equation}

\noindent From ~\eqref{soleq2}, $\forall \tau \in \mathcal{R}$, if
$\lambda_{\tau}(z)=0$ the polynomial
$\left(\sum_{i=1}^n(-1)^{i+j}M_{i,j}(\partial \tilde f/
\partial z_j)\right)_z(Z)$ is divisible by
$\displaystyle{ \lambda_{\tau}(Z)} $, and since the $\lambda_{\tau}$
are pairwise relatively prime, $\forall \mathcal{U}\subseteq
\mathcal{R}$ if $z\in \bigcap_{\tau\in \mathcal{U}} H_{\tau}$, then
$\left(\sum_{i=1}^n(-1)^{i+j}M_{i,j}(\partial \tilde f/
\partial z_j)\right)_z(Z)$ is divisible by
$ \prod_{\tau \in \mathcal{U}}\lambda_{\tau}(Z)$.\vskip 5pt

 In the reducible case, in convenient bases, the jacobian matrix of $P$ is block
 diagonal. We have $J_P=J_{P^1}\ldots J_{P^s}$ and
 $M_{i,j}$is of the form $J_{P^1}\ldots J_{P^{l-1}}M_{i,j}^lJ_{P^{l+1}}\ldots J_{P^s}$
 where $M_{i,j}^l$
 is a minor of the $n_l\times n_l$ block associated with $P^l$, Chevalley map of $W^l$.
 After simplification we see that it is sufficient to study each block.

For any finite reflection group: \[h=k_n=\sum_{1\le j\le
n}(k_j-1)-\sum_{1\le j\le n-1}(k_j-1)+1,\] where the first sum
$\sum_{1\le j\le n}(k_j-1)$ is the degree of the jacobian
determinant, equal to the number of linear forms $\lambda_{\tau}$
which is the number $d={\mathcal R}^{\#}$ of reflections in $W$. The
second sum $\sum_{1\le j\le n-1}(k_j-1)$ is the least degree $s$ of
the minors $M_{i,j}$ of the jacobian determinant of system
~\eqref{sys1}.

If $W$ is reducible, the formula also holds for each irreducible
component. In particular for any component with the greatest Coxeter
number $h$, we have $h=1+d'-s'$ where $d'$ and $s'$ are those
associated with this irreducible component. Accordingly, we may and
will assume from now on, without loss of generality, that $W$ is an
irreducible Coxeter group.

\subsection*{Stratification and Compensation by the $M_{i,j}$.}

 We denote by $S_p$ a stratum in $P^{-1}({\bf R}^n)$, which is a connected component
 of the intersection of $\Gamma$ and exactly $p$ of the reflecting hyperplanes.
 The points of each stratum are stabilized by the same isotropy group,
 subgroup of $W$ generated by reflections about the hyperplanes containing the
 stratum.
 The different possible isotropy subgroups and strata
 types may be determined from the Dynkin diagram.  The stratum of dimension
 $0$ is the origin. The strata of dimension 1 are those determined by
 removing only one point in the Dynkin diagram, they are strata $S_q$ such
 that their closure is $\overline{S_q}=S_q\cup\{0\}$. At the other end the
 strata of dimension $n$ are the connected components of the regular set
 in $\Gamma$.  \medskip

  On any stratum $S_p$, by lemma \ref{lem1}
  and lemma \ref{lem2}, $({\displaystyle{\partial \widetilde
  F}/ \displaystyle{\partial p_j}})\circ P \in \mathcal{H}^{r-1+s_p-p}(S_p)$,
  if $M_{i,j}$ is at least $(s_p-1)$-flat on $S_p$.
  In a neighborhood of $z\in S_p$ we have $P=q\circ v$
 with $q$ invertible and $v$ the Chevalley mapping of $W_z$,
 isotropy subgroup of $z$ (and of any point in $S_p$). Observe that $W_z$ is
 not irreducible: it is not essential since it stabilizes $S_p$, and
 it may also have several irreducible components.
 The adjoint of the jacobian matrix of $P$ which is the transpose of its
 comatrix is the product in this order of the adjoint of the jacobian
 matrices $V$ and $Q$ of $v$ and $q$ respectively.
 We have $M_{i,j}=\sum_{k=1}^n V_{k,j}Q_{i,k}$. The $V_{k,j}$ and accordingly
 the $M_{i,j}$ are $(s_z-1)$-flat on $S_p$. So $s_p=s_z$ and,
 since $p$ is the number $d_z$ of reflections in
 $W_z,\; 1-s_p+p$ is the Coxeter number $h_z=1-s_z+d_z$ of $W_z$.

 The isotropy group of the points $z\in S_p$ is a subgroup of
 the isotropy group of the points $z'\in S_{p+q}$. Therefore $h_z\le
 h_{z'}$, or $1+p-s_p\le 1+p+q-s_{p+q}$.
 Also $1+d-s_d=h$ is larger than
 $h_z=1+d_z-s_z$ for any $z\neq 0$.

 Finally, the minors $M_{i,j}$ are homogeneous polynomials of
degree: \[s_j=\displaystyle{\sum_{1\le u\le n, u\neq j} (k_u-1)} \ge
  s=\displaystyle{ \sum_{1\le u\le n-1} (k_u-1)}.\] They are at least
  $(s-1)$-flat on the intersection of the $\mathcal{R}^{\#}=d$
 reflecting hyperplanes. (Observe that this might be used to get
 $s_p=s_z$ by induction).

 In ~\eqref{soleq2}, lemma \ref{lem4} applies to the closure
 of each connected component of the regular set, with $\displaystyle{A_i=
{\partial \tilde f / \partial z_i} \in \mathcal{H}^{r-1}(P^{-1}({\bf
R}^n))}$, $Q_i= (-1)^{i+j}M_{i,j}$,
 and gives the $[(r-1)+s-d]$-continuity of the ${\displaystyle{(\partial
\widetilde F}/ \displaystyle{\partial p_j})}\circ P$ on their union
$P^{-1}({\bf R}^n)$.
 The result we needed to complete the proof:
 \[{\displaystyle{\partial \widetilde F}\over \displaystyle{\partial
 p_j}}\circ P \in \mathcal{H}^{r-1-d+s}(P^{-1}({\bf R}^n)),\qquad 1+d-s=h\]
 is now a consequence of lemma \ref{lem3}.$ \footnote{The closure of each connected
 component of the regular set being convex and thus Whitney
 1-regular (\cite{24}), lemma \ref{lem4} could directly give
 the result.}.\qquad \square$

 \begin{rem} The above result gives a loss
 of differentiability of $1+d-s$, where $s$ is the least degree
 of the $M_{i,j}$. Actually $M_{i,j}$ is the jacobian of the polynomial
 mapping:
 \[(z_1, \ldots , z_{i-1}, z_{i+1},\ldots , z_n;\; z_i)\mapsto
 (p_1(z),\ldots, p_{j-1}(z), p_{j+1}(z),\ldots,p_n(z);\; z_i).\]

\noindent This mapping is invariant by the subgroup $W_i$ of $W$
that leaves invariant the $i^{th}$ coordinate axis in ${\bf R}^n$,
say ${\bf R\; e}_i$ (\cite{3}). This subgroup $W_i$ is generated by
the subset $\mathcal{R}_i\subset \mathcal{R} $ of the reflections it
contains. These are the reflections $\alpha$ in $W$ such that
$\alpha({\bf e}_i)={\bf e}_i$, about the hyperplanes $H_{\alpha}$
containing ${\bf e}_i$ \footnote{The description of $W_i$ given in
\cite{3} was correct. Unfortunately in \cite{4} it was not. Although
not essential in the reasoning it was misleading. The explicit
computations were correct however and gave the best result in the
cases of $A_n$, $B_n$ and $I_2(k)$.}. The $M_{i,j},\; j=1,\ldots,
l,$ as jacobians of $W_i$-invariant polynomial mappings are
polynomial multiples of $(\prod_{\tau\in \mathcal{R}_i}
\lambda_{\tau})$. In \cite{4} the formula for the loss of
differentiability at each step was also of the form $1+d-s$, but
 $s$ was $\min_{1\le i\le n}{\mathcal R}_i^{\#}$. Clearly
 $\min_{1\le i\le n}{\mathcal R}_i^{\#}\le \min_{1\le i,j\le n}
 {{\rm degree}M_{i,j}}$. In some cases ($A_n, B_n, I_2(k)$)
 the equality holds but in general the loss of differentiability
 given by \cite{4} was not the best possible. Considering $H_3$ for an example,
 we now have $s= 6$ instead of $2$, and
 the class of differentiability of $F$ is $[r/10]$ instead of $[r/14]$.
\end{rem}

 All the operations from $f\in \mathcal{C}^r({\bf
R}^n)^W$ up to $F\in \mathcal{C}^{[r/h]}({\bf R}^n)$ are linear and
continuous when using the natural Fr\'echet topologies. A modulus of
continuity for the Whitney conditions could be followed from
$\|f\|^r$ to $\|F\|^{[r/h]}$. So Chevalley's theorem in class
$\mathcal{C}^r $ may be restated as:

{\bf Theorem 1.1.} {\it Let $W$ be a finite group generated by
reflections acting orthogonally on ${\bf R}^n$, $P$ the Chevalley
polynomial mapping associated with $W$, and $h=k_n$ the highest
degree of the coordinate polynomials in $P$ (equal to the greatest
Coxeter number of the irreducible components of $W$). There exists a
linear and continuous mapping:
$$\mathcal{C}^r({\bf R}^n)^W \ni f \to F\in \mathcal{C}^{[r/h]}
({\bf R}^n)$$ such that $f=F\circ P$.}

\begin{rem} Theorem \ref{thm1} gives global results. As it is clear
from the proof, in the neighborhood of each point $x$ the loss of
differentiability is governed by the isotropy group $W_x$ and its
Coxeter number $h_x$.

About the partial derivatives at the origin (or a point $x$ where
$W_x=W$), since the $M_{i,j}$ are homogeneous of degree $s_j=
\sum_{1\le u\le n, u\neq j}(k_u-1)$, we see that $({\partial
\widetilde F/\partial p_j})\circ P$ is of class ${\mathcal
H}^{r-k_j}$. Reasoning as above, we could show that the partial
derivatives $\displaystyle{{\partial^{|m|} \widetilde F/
\partial P^m}}$ of order $m=(m_1,\ldots,m_n)$ are continuous if
$m_1k_1+\ldots +m_nk_n\le r$. For instance the partial derivatives
in $W$-invariant directions are continuous up to the order $r$.
\end{rem}

\section{Counter Example.}

 Let us give a counter example which applies to almost every finite
 reflection group. It is sufficient to consider essential
 irreducible groups.

We consider $F:{\bf R}^n \to {\bf R}$ defined by $
F(y)=y_n^{s+\alpha}$ for some integer $s$ and an $\alpha\in ]0,1[$.
$F$ is of class $\mathcal{C}^s$ but not of class $\mathcal{C}^{s+1}$
in any neighborhood of $\{y | y_n = 0\}$. Let $P$ be the Chevalley
mapping associated with some finite irreducible Coxeter group $W$
acting on ${\bf R}^n$ and consider the composite mapping $F\circ
P(x)= p_n^{s+\alpha}(x)$. We study the differentiability of this
mapping when $p_n(x)=0$.

A set of basic invariants is available in \cite{18} for any finite
Coxeter group. Disregarding $D_n$, for any other group there exists
an invariant set of linear forms $\{L_1,\ldots ,L_v\}$ the kernels
of which intersect only at the origin, and such that for
$i=1,\ldots,n,\;p_i(X)=\sum_{j=1}^v [L_j(X)]^{k_i}$ with $k_i$s as
determined in [7]. With the two exceptions of $A_{2n}$ and
$I_2(2p+1)$, $k_n$ is even and therefore $p_n(x)$ vanishes only at
the origin. We will not study the two exceptional cases, but a
fairly general counter example is given in [1] for symmetric
functions and thus for $A_n$ (including $A_2=I_2(3)$). As usual,
$D_n$ does not follow the general line but we may choose $p_n(x)
=\sum_1^n x_i^{2(n-1)}$ and the results of the general case apply.

We have $p_n(x)=\sum_1^v[L_i(x)]^{k_n}$, and since $|L_i(x)|\le
a_i|x|,\; i=1,\ldots , v$ for some numerical constants $a_i$, we
have the estimate $|p_n(x)|\le (\sum_1^v
a_i^{k_n})|x|^{k_n}=A|x|^{k_n}$.

Analogously, since $|D^1L_i(x)|\le b_i$ for some numerical constants
$b_i$, we get: $$|D^jp_n(x)|\le \sum_1^v b_i^j{k_n\choose
j}|L_i(x)|^{k_n-j}= B_j|x|^{k_n-j}.$$

The derivatives of the composite mapping $p_n^{s+\alpha}(x)$ are
given by the Faa di Bruno formula: \[D^k p_n^{s+\alpha}(x)= \sum
{k!\over \mu_1!\ldots \mu_q!}
D^py_n^{s+\alpha}(p_n(x))\big({D^1p_n(x)\over 1!}\big)^{\mu_1}
\ldots \big({D^qp_n(x)\over q!}\big)^{\mu_q},\] where the sum is
over all the $q$-tuples $(\mu_1,\ldots\mu_q)\in {\bf N}^q$ such that
$1\mu_1+\ldots +q\mu_q=k$, with $p=\mu_1+\ldots +\mu_q$. There are
constants $C_{(\mu_1,\ldots,\mu_q)}$ such that:
$$ |\big({D^1p_n(x)\over 1!}\big)^{\mu_1} \ldots
\big({D^qp_n(x)\over q!}\big)^{\mu_q}|\le
C_{(\mu_1,\ldots,\mu_q)}|x|^{(k_n-1)\mu_1+\ldots
+(k_n-q)\mu_q}=C_{(\mu_1,\ldots,\mu_q)}|x|^{k_np-k},$$ and therefore
constants $A_{(\mu_1,\ldots\mu_q)}$ and $A$ such that: \[|D^k
p_n^{s+\alpha}(x)|\le \sum
A_{(\mu_1,\ldots\mu_q)}|x|^{k_n(s+\alpha-p)}|x|^{k_np-k} \le A
|x|^{k_ns+k_n\alpha-k}.\] This shows that the derivatives of order
$k\le k_ns$ tend to $0$ at the origin while the derivatives of order
$k_ns+1$ will not if $\alpha<1/k_n$. This means that the composite
mapping $f=F\circ P$ is of class $\mathcal{C}^{k_n s}$ but not of
class $\mathcal{C}^{k_n s+1}$ at $x=0$ and it factors through $F$
which is of class $\mathcal{C}^s$ and not of class
$\mathcal{C}^{s+1}$. The loss of differentiability is as given in
theorem \ref{thm1} and cannot be reduced.

% ----------------------------------------------------------------
\bibliographystyle{amsplain}

\end{document}